\documentclass[a4paper,11pt]{article}

\usepackage[english]{babel}
\usepackage{amsmath}
\usepackage{amssymb}
\usepackage{fancyhdr}
\usepackage[latin1]{inputenc}
\usepackage[T1]{fontenc}
\usepackage{stmaryrd}
\usepackage{amsthm}
\usepackage{graphicx}
\usepackage{geometry}

\newtheorem{thrm}{Theorem}[subsection]
\newtheorem{prpstn}{Proposition}[subsection]
\newtheorem{lmm}{Lemma}[subsection]
\newcommand{\e}{\varepsilon}
\newcommand{\E}{\mathbb{E}}
\newcommand{\dd}{\mathrm{d}}
\newcommand{\OO}{\mathrm{O}}

\newcommand{\ii}{\mathrm{i}}

\makeatletter
\renewcommand{\thefigure}{\ifnum \c@section>\z@ \thesection.\fi
 \@arabic\c@figure}
\@addtoreset{figure}{section}
\makeatother

\newcommand{\auteur}[4]{ \author{ { #1}\footnote{#2} \vspace{0.5cm} 
 \\  
 \hspace{-2.0cm} \small \noindent\begin{tabular}{l}   #3 \\ #4\end{tabular}}}

\auteur{Fanny Godet}{fanny.godet@math.univ-nantes.fr \endgraf I would like to thank Anne Philippe, my PhD advisor, who monitored my work.}{Laboratoire de Mathématiques Jean Leray, UMR CNRS 6629  } { Université de Nantes
2 rue de la Houssinière - BP 92208 
F-44322 Nantes Cedex 3 }
\begin{document}
\title{Linear Prediction of Long-Memory Processes:\\
Asymptotic Results on Mean-squared Errors}
\maketitle
\begin{abstract}
We present two approaches for linear prediction of long-memory time series. The first approach consists in truncating the Wiener-Kolmogorov predictor by restricting the observations to the last $k$ terms, which are the only available values in practice. We derive the asymptotic behaviour of the mean-squared error as $k$ tends to $ + \infty$. By contrast, the second approach is non-parametric. An AR($k$) model is fitted to the long-memory time series and we study the error that arises in this misspecified model.
\end{abstract}


\par ARMA (autoregressive moving-average) processes are often called short-memory processes because their covariances decay rapidly (i.e.\ exponentially). On the other hand, a long-memory process is characterised by the following feature: the autocovariance function $\sigma$  decays more slowly i.e.\ it is not absolutely summable. They are so-named because of the strong association between observations widely separated in time. The long-memory time series models have attracted much attention lately and there is now a growing realisation that time series possessing long-memory characteristics arise in subject areas as diverse as Economics, Geophysics, Hydrology or telecom traffic (see, e.g., \cite{mandelbrot} and \cite{granger}). Although there exists substantial literature on the prediction of short-memory processes(see \cite{bhansali} for the univariate case or \cite{lewis} for the multivariate case), there are fewer results for long-memory time series. In this paper, we consider the question of the prediction of the latter.
\par More precisely, we compare two prediction methods for long-memory process. Our goal is a linear predictor of $X_{k+h}$ based on observed time points which is optimal in the sense that it minimizes the mean-squared error. The paper is organized as follows. First we introduce our model and our main assumptions. Then, in section \ref{meilleurpredicteurlineairetronque}, we study the best linear predictor i.e.\ the Wiener-Kolmogorov predictor proposed by \cite{whittle} and by \cite{tronc} for long-memory time series. In practice, only the last $k$ values of the process are available. Therefore we need to truncate the infinite series in the definition of the predictor and to derive the asymptotic behaviour of the mean-squared error as $k$ tends to $ + \infty$.
\par In Section \ref{sectionarp} we discuss the asymptotic properties of the forecast error if we fit a misspecified AR($k$) model to a long-memory time series. This approach has been proposed by \cite{arp} for fractional noise series F($d$). His simulations show that high-order AR-models predict fractional integrated noise very well. 
\par Finally in Section \ref{prevh} we compare the two previous methods for $h$-step prediction. We give some asymptotic properties of the mean-squared error of the linear least-squares predictor as $h$ tends to $ + \infty$ in the particular case of long-memory processes. Then we study our $k$-th order predictors order as $k$ tends to $ + \infty$.

\section{Model}
\label{intro}
\par Let $(X_n)_{n \in \mathbb{Z}}$ be a discrete-time (weakly) stationary process in $\mathrm{L}^2$ with mean 0 and $\sigma$ its autocovariance function. We assume that the process $(X_n)_{n \in \mathbb{Z}}$ is a long-memory process i.e.:
\begin{displaymath}
\sum_{k=-\infty}^{\infty} |\sigma(k)|=\infty.
\end{displaymath}
The process $(X_n)_{n \in \mathbb{Z}}$ admits an infinite moving average representation as follows:
\begin{equation}
 X_n=\sum_{j=0}^{\infty} b_j \e_{n-j} \label{defbj}
\end{equation}
where $(\varepsilon_n)_{n\in \mathbb{Z}}$ is a white-noise series consisting of uncorrelated random variables, each with mean zero and variance $\sigma^2_{\varepsilon}$ and $(b_j)_{j \in \mathbb{N}}$ are square-summable.
 We shall further assume that $(X_n)_{n \in \mathbb{Z}}$ admits an infinite autoregressive representation:
\begin{equation}
\e_n=\sum_{j=0}^{\infty} a_j X_{n-j}, \label{defaj}
\end{equation}
where the $(a_j)_{j \in \mathbb{N}}$ are absolutely summable.
We assume also that $(a_j)_{j \in \mathbb{N}}$ and  $(b_j)_{j \in \mathbb{N}}$, occurring respectively in \eqref{defaj} and \eqref{defbj}, satisfy the following conditions for all $\delta>0$:
\begin{eqnarray}
|a_j|&\leq&C_1  j^{-d-1+\delta} \label{aj} \\
|b_j|&\leq& C_2 j^{d-1+\delta} .\label{bj}
\end{eqnarray}
where $C_1$ and $C_2$ are constants and $d$ is a parameter verifying $d \in ]0,1/2[$. For example, a FARIMA process $(X_n)_{n \in \mathbb{Z}}$ is the stationary solution to the difference equations:
\begin{displaymath}
\phi(B) (1-B)^d X_n = \theta(B) \e_n
\end{displaymath}
where $(\e_n)_{n \in \mathbb{Z}}$ is a white noise series, $B$ is the backward shift operator and $\phi$ and $\theta$ are polynomials with no zeroes on the unit disk. Its coefficients verify 
\begin{eqnarray}
|a_j|&\underset{+\infty}\sim& C_1 j^{-d-1}  \nonumber \\
|b_j|&\underset{+\infty}\sim& C_2 j^{d-1} \nonumber 
\end{eqnarray}
and thus \eqref{aj} and \eqref{bj} hold. When $\phi=\theta=1$, the process $(X_n)_{n \in \mathbb{Z}}$ is called \textit{fractionally integrated noise} and denoted F($d$). More generally, $(a_j)_{j \in \mathbb{N}}$ and  $(b_j)_{j \in \mathbb{N}}$ verify conditions \eqref{aj} and \eqref {bj} if:
\begin{eqnarray}
|a_j|&\underset{+\infty}\sim& L(j) j^{-d-1}  \nonumber \\
|b_j|&\underset{+\infty}\sim& L'(j) j^{d-1} \nonumber 
\end{eqnarray}
where $L$ and $L'$ are slowly varying functions.
A positive function $L$ is a slowly varying function in the sense of \cite{zygmund} if, for any $\delta>0$, $x \mapsto x^{-\delta}L(x)$ is decreasing and $x \mapsto x^{\delta}L(x)$ is increasing. 
\par  The condition \eqref{bj} implies that the autocovariance function $\sigma$ of the process $(X_n)_{n \in \mathbb{Z}}$ verifies:
\begin{equation}
\forall \delta>0, \exists C_3 \in \mathbb{R} ,\quad \vert \sigma(j)\vert \leq  C_3j^{2d-1+\delta} .\label{sigma}
\end{equation}
Notice that it suffices to prove \eqref{sigma} for $\delta$ near 0 in order to verify  \eqref{sigma} for $\delta >0$ arbitrarily chosen. So we prove \eqref{sigma} for $\delta<\frac{1-2d}{2}$:
\begin{eqnarray*}
\sigma(k)&=& \sum_{j=0}^{+\infty} b_j b_{j+k} \\
\vert \sigma(k) \vert &\leq&  \sum_{j=1}^{+\infty}\vert b_j b_{j+k}\vert+\vert b_0b_k \vert \\
&\leq& C_2^2\sum_{j=1}^{+\infty} j^{d-1+\delta} (k+j)^{d-1+\delta} +\vert b_0b_k \vert\\
&\leq& C_2^2 \int_{0}^{+\infty} j^{d-1+\delta} (k+j)^{d-1+\delta} \dd j +\vert b_0b_k \vert\\
&\leq& C_2^2k^{2d-1+2\delta}\int_{0}^{+\infty}j^{d-1+\delta} (1+j)^{d-1+\delta} \dd j+C_2k^{d-1+\delta}\\
&\leq& C_3k^{2d-1+2\delta}
\end{eqnarray*}
 More accurately, \cite{inoueregularly} has proved than if:
\begin{displaymath}
b_j \sim L\left( j\right) j^{d-1} 
\end{displaymath}
then 
\begin{displaymath}
\sigma(j) \sim j^{2d-1} \left[ L\left( j\right)\right]^2 \beta(1-2d,d)
\end{displaymath}
where $L$ is a slowly varying function and $\beta$ is the beta function. The converse is not true, we must have more assumptions about the series $(b_j)_{j\in \mathbb{N}}$ in order to get an asymptotic equivalent for $(\sigma(j))_{j\in \mathbb{N}}$ (see \cite{inouepartialautocorr}).

\section{Wiener-Kolmogorov Next Step Prediction Theory}
\label{meilleurpredicteurlineairetronque}
\subsection{Wiener-Kolmogorov Predictor}
\par The aim of this part is to compute the best linear one-step predictor (with minimum mean-square distance from the true random variable) knowing all the past $\{X_{k+1-j},j \leqslant 1\}$. Our predictor is therefore an infinite linear combination of the infinite past:
\begin{equation}
\widetilde{X_k}(1)=\sum_{j=0}^{\infty} \lambda(j) X_{k-j} \nonumber 
\end{equation}
 where $(\lambda(j))_{j \in \mathbb{N}}$ are chosen to ensure that the mean squared prediction error:
\begin{displaymath}
\mathbb{E}\big[ \big(\widetilde{X_k}(1)-X_{k+1}\big)^2\big]
\end{displaymath}
is as small as possible. Following \cite{whittle}, and in view of the moving average representation of $(X_n)_{n \in \mathbb{Z}}$, we may rewrite our predictor $\widetilde{X_k}(1)$ as:
\begin{displaymath}
\widetilde{X_k}(1)=\sum_{j=0}^{\infty}\phi(j)\varepsilon_{k-j}.
\end{displaymath}
where $(\phi(j))_{j \in \mathbb{N}}$ depends only on $(\lambda(j))_{j \in \mathbb{N}}$ and $(\epsilon_n)_{n \in \mathbb{Z}}$ and $(a_j)_{j \in \mathbb{N}}$ are defined in \eqref{defaj}.
From the infinite moving average representation of $(X_n)_{n \in \mathbb{Z}}$ given below in \eqref{defbj}, we can rewrite the mean-squared prediction error as:
\begin{eqnarray*}
\mathbb{E}\big[ \big(\widetilde{X_k}(1)-X_{k+1}\big) ^2\big]
&=&\mathbb{E}\left[ \left( \sum_{j=0}^{\infty}\phi(j)\varepsilon_{k-j}-\sum_{j=0}^{\infty} b(j)\varepsilon_{k+1-j}\right) ^2\right]  \\  
&=&\mathbb{E}\left[ \left(\varepsilon_{k+1}-\sum_{j=0}^{\infty}\left( \phi(j)-b(j+1)\right) \varepsilon_{k-j}\right) ^2\right] \\
&=&\left( 1+\sum_{j=0}^{\infty}\big( b_{j+1}-\phi(j)\big) ^2\right) \sigma_{\e} ^2
\end{eqnarray*}
since the random variables $(\e_n)_{n \in \mathbb{Z}}$ are uncorrelated with variance $\sigma_{\e}^2$. The smallest mean-squared prediction error is obtained when setting $ \phi(j)=b_{j+1}$ for $j\geq 0$. 
 \par The smallest prediction error of $(X_n)_{n \in \mathbb{Z}}$ is $\sigma_{\e}^2$ within the class of linear predictors.
Furthermore, if 
\begin{displaymath}
A(z)=\sum_{j=0}^{+\infty} a_j z^j,
\end{displaymath}
denotes the characteristic polynomial of the  $(a(j))_{j \in \mathbb{Z}}$ and
\begin{displaymath}
B(z)=\sum_{j=0}^{+\infty} b_j z^j,
\end{displaymath}
that of the $(a(j))_{j \in \mathbb{Z}}$, then in view of the identity, $A(z)=B(z)^{-1},\,\vert z \vert \leq 1$, we may write:
\begin{equation}
\widetilde{X_k}(1)=-\sum_{j=1}^{\infty}a_jX_{k+1-j} .\label{serieinf}
\end{equation}

\subsection{Mean Squared Prediction Error when the Predictor is Truncated}
\label{erreurtronc}
\par In practice, we only know a finite subset of the past, the one which we have observed. So the predictor should only depend on the observations. Assume that we only know the set $\{X_1, \ldots, X_k\}$ and that we replace the unknown values by 0, then we have the following new predictor: 
\begin{equation}
\widetilde{X'_k}(1)=-\sum_{j=1}^{k}a_jX_{k+1-j} .\label{pred1}
\end{equation}
It is equivalent to say that we have truncated the infinite series \eqref{serieinf} to $k$ terms.
The following proposition provides us the asymptotic properties of the mean squared prediction error as a function of $k$.

\begin{prpstn} \label{erreurtronque} Let $(X_n)_{n \in \mathbb{Z}}$ be a linear stationary process defined by \eqref{defbj}, \eqref{defaj} and verifying conditions \eqref{aj} and \eqref{bj}. We can approximate the mean-squared prediction error of $\widetilde{X'_k}(1)$ by:
\begin{displaymath}
\forall \delta>0, \quad \mathbb{E}\big( \big[X_{k+1}-\widetilde{X'_k}(1)\big]^2 \big)  =\sigma_{\e}^2+\OO(k^{-1+\delta}).
\end{displaymath}
Furthermore, this rate of convergence $\OO(k^{-1})$ is optimal since for fractionally integrated noise, we have the following asymptotic equivalent:
\begin{displaymath}
\mathbb{E}\big( \big[X_{k+1}-\widetilde{X'_k}(1)\big]^2 \big)  =\sigma_{\e}^2+C k^{-1}+\mathrm{o}\left(k^{-1} \right) .
\end{displaymath}

\end{prpstn}
Note that the prediction error is the sum of $\sigma_{\e}^2$, the error of Wiener-Kolmogorov model and the error due to the truncation to $k$ terms which is bounded by $\OO(k^{-1+\delta})$ for all $\delta>0$.
\begin{proof}
\begin{eqnarray}
X_{k+1}-\widetilde{X'_k}(1) &=& X_{k+1}-\widetilde{X_k}(1)+\widetilde{X_k}(1)-\widetilde{X'_k}(1) \nonumber \\
&=& X_{k+1} -\sum_{j=0}^{+\infty}b_{j+1} \e_{k-j} -\sum_{j=k+1}^{+\infty}a_j X_{k+1-j} \nonumber \\
&=& \e_{k+1}-\sum_{j=k+1}^{+\infty}a_j X_{k+1-j}  \label{orthogonalite}.
\end{eqnarray}
The two parts of the sum \eqref{orthogonalite} are orthogonal for the inner product associated with the mean square norm. Consequently:
\begin{displaymath}
\mathbb{E}\big( \big[X_{k+1}-\widetilde{X'_k}(1)\big]^2 \big) 
 = \sigma_{\e}^2+ \sum_{j=k+1}^{\infty} \sum_{l=k+1}^{\infty} a_j a_l \sigma(l-j)  .
\end{displaymath}
For the second term of the sum we have:
\begin{eqnarray}
\bigg \vert \sum_{j=k+1}^{+\infty} \sum_{l=k+1}^{+\infty} a_j a_l \sigma(l-j) \bigg \vert&=&\bigg \vert 2 \sum_{j=k+1}^{+\infty} a_j \sum_{l=j+1}^{+\infty}a_l \sigma(l-j) + \sum_{j=k+1}^{+\infty}a_j^2 \sigma(0)\bigg \vert \nonumber \\
&\leq & 2 \sum_{j=k+1}^{+\infty}|a_j|\left|a_{j+1}\right|| \sigma(1)|+ \sum_{j=k+1}^{+\infty}a_j^2 \sigma(0) \nonumber\\
&& +2\sum_{j=k+1}^{+\infty} |a_j|\sum_{l=j+2}^{+\infty}|a_l|| \sigma(l-j)|\nonumber
\end{eqnarray}
from the triangle inequality, it follows that:
\begin{eqnarray}
&&\bigg \vert \sum_{j=k+1}^{+\infty} \sum_{l=k+1}^{+\infty} a_j a_l \sigma(l-j) \bigg \vert \nonumber \\
&\leq &
C_1^2C_3\left( 2\sum_{j=k+1}^{+\infty}j^{-d-1+\delta}(j+1)^{-d-1+\delta}+\sum_{j=k+1}^{+\infty}\left( j^{-d-1+\delta}\right)^2 \right) \label{sommesimple}\\
&+&2C_1^2C_3\sum_{j=k+1}^{+\infty}j^{-d-1+\delta} \sum_{l=j+2}^{+\infty}l^{-d-1+\delta}\vert l-j\vert^{2d-1+\delta} \label{sommedouble} 
\end{eqnarray}
 for all $\delta>0$ from inequalities \eqref{aj} and \eqref{sigma}. Assume now that $\delta <1/2-d$. For the terms \eqref{sommesimple}, since $j \mapsto j^{-d-1+\delta} (j+1)^{-d-1+\delta}$ is a positive and decreasing function on $\mathbb{R}^+$, we have the following approximations:
 \begin{eqnarray*}
 2C_1^2C_3\sum_{j=k+1}^{+\infty}j^{-d-1+\delta}(j+1)^{-d-1+\delta} & \sim &2C_1^2C_3 \int_k^{+\infty}j^{-d-1+\delta}(j+1)^{-d-1+\delta} \dd j \\
  & \sim &\frac{2C_1^2C_3}{1+2d-2\delta} k^{-2d-1+2\delta} 
 \end{eqnarray*}
Since the function $j \mapsto \left( j^{-d-1+\delta}\right) ^2$ is also positive and decreasing, we can establish in a similar way that:
 \begin{eqnarray*}
 C_1^2C_3\sum_{j=k+1}^{+\infty}\left( j^{-d-1+\delta}\right) ^2 & \sim&C_1^2C_3\int_k^{+\infty}\left( j^{-d-1+\delta}\right) ^2\dd j \\
 & \sim &\frac{C_1^2C_3}{1+2d-2\delta} k^{-2d-1+2\delta}. 
 \end{eqnarray*}

For the infinite double series \eqref{sommedouble}, we will similarly compare the series with an integral. In the next Lemma, we establish the necessary result for this comparison:
\begin{lmm} \label{sommedoublecomp}
Let $g$ the function $(l,j) \mapsto j^{-d-1+\delta} \, l^{-d-1+\delta}\,\vert l-j\vert^{2d-1+\delta}$. Let $m$ and $n$ be two positive integers. We assume that $\delta<1-2d $ and $m \geq \frac{\delta -d-1}{\delta+2d-1}$ for all $\delta \in \left] 0,\frac{\delta -d-1}{\delta+2d-1}\right[ $. 
We will call $A_{n,m}$ the square $[n,n+1]\times[m,m+1]$. If $n \geq m+1$ then
\begin{displaymath}
\int_{A_{n,m}}g(l,j)\, \dd j \,\dd l \geq g(n+1,m).
\end{displaymath}
\end{lmm}

Assume now that $\delta<1-2d $ without loss of generality.
Thanks to the previous Lemma and the asymptotic equivalents of \eqref{sommesimple}, there exists $K \in \mathbb{N}$ such that if $k>K$:
\begin{eqnarray}
\bigg \vert \sum_{j=k+1}^{+\infty} \sum_{l=k+1}^{+\infty} a_j a_l \sigma(l-j) \bigg \vert&\leq & C \int_{k+1}^{+\infty} j^{-d-1+\delta}\left[  \int_{j}^{+\infty} l^{-d-1+\delta} (l-j)^{2d-1+\delta} \dd l \right]  \dd j +\OO\left( k ^{-2d-1+2\delta}\right) \nonumber
\end{eqnarray}
By using the substitution $jl'=l$ in the integral over $l$ we obtain:
\begin{displaymath}
\bigg \vert \sum_{j=k+1}^{+\infty} \sum_{l=k+1}^{+\infty} a_j a_l \sigma(l-j) \bigg \vert\leq  C'\int_{k+1}^{+\infty}j^{-2+3\delta} \int_{1}^{+\infty} l^{-d-1+\delta} (l-1)^{2d-1+\delta} \dd l \dd j +\OO\left( k ^{-2d-1}\right).  
\end{displaymath} Since if $\delta<(1-d)/2$
\begin{displaymath}
\int_1^{+\infty} l^{-d-1+\delta}(l-1)^{2d-1+\delta} \dd l < +\infty,
\end{displaymath}
it follows:
\begin{eqnarray}
\bigg \vert \sum_{j=k+1}^{+\infty} \sum_{l=k+1}^{+\infty} a_j a_l \sigma(l-j) \bigg \vert&\leq & \OO\left(k^{-1+3\delta} \right)+\OO\left( k ^{-2d-1}\right) \nonumber\\
&\leq & \OO\left(k^{-1+3\delta} \right).
\end{eqnarray}

If $\delta>0$, $\delta<1-2d $ and $\delta<(1-d)/2$, we have:
\begin{displaymath}
\bigg \vert \sum_{j=k+1}^{+\infty} \sum_{l=k+1}^{+\infty} a_j a_l \sigma(l-j) \bigg \vert = \OO\left(k^{-1+3\delta} \right).
\end{displaymath}
Notice that if the equality is true under the assumptions $\delta>0$, $\delta<1-2d $ and $\delta<(1-d)/2$, it is also true for any $\delta>0$. Therefore we have proven the first part of the theorem. \\
 We prove now that there exists long-memory processes whose prediction error attains the rate of convergence $k^{-1}$. Assume now that $(X_n)_{n \in \mathbb{Z}}$ is fractionally integrated noise F($d$), which is the stationary solution of the difference equation:
\begin{equation}
X_n=(1-B)^{-d} \e_n \label{Fd}
\end{equation}
with $B$ the usual backward shift operator, $(\e_n)_{n \in \mathbb{Z}}$ is a white-noise series and $d\in \left]0, 1/2\right[ $ (see for example \cite{bd}). We can compute the coefficients and obtain that:
\begin{displaymath}
\forall j>0,\quad a_j=\frac{\Gamma(j-d)}{\Gamma(j+1)\Gamma(-d)} \textrm{ and } \forall j\geq 0,\quad \sigma(j)=\frac{(-1)^j\Gamma(1-2d)}{\Gamma(j-d+1)\Gamma(1-j-d)} \sigma_\e^2
\end{displaymath}
then we have:
\begin{displaymath}
\forall j>0,\quad a_j<0\textrm{ and } \forall j\geq 0,\quad \sigma(j)>0
\end{displaymath}
and 
\begin{displaymath}
a_j \sim \frac{j^{-d-1}}{\Gamma(-d)}\textrm{ and }\sigma(j)\sim \frac{j^{2d-1}\Gamma(1-2d)}{\Gamma(d)\Gamma(1-d)} \quad \textrm{when } j \rightarrow \infty. 
\end{displaymath}

In this particular case, we can estimate the prediction error more precisely:
\begin{eqnarray}
\sum_{k+1}^{+\infty} \sum_{k+1}^{+\infty} a_j a_l \sigma(l-j)&=&\sum_{k+1}^{+\infty} |a_j|\sum_{j+1}^{+\infty}|a_l|| \sigma(l-j)| + \sum_{k+1}^{+\infty}a_j^2 \sigma(0) \nonumber\\
&\sim&\frac{\Gamma(1-2d)}{\Gamma(-d)^2\Gamma(d)\Gamma(1-d)}\int_{k+1}^{+\infty}j^{-2} \int_{1/j+1}^{+\infty} l^{-d-1} (l-1)^{2d-1} \dd l \dd j +\OO\left( k ^{-2d-1}\right) \nonumber\\
\sum_{k+1}^{+\infty} \sum_{k+1}^{+\infty} a_j a_l \sigma(l-j) &\sim& \frac{\Gamma(1-2d)\Gamma(2d)}{\Gamma(-d)^2\Gamma(d)\Gamma(1+d)} k^{-1} \label{premest}
\end{eqnarray}
The asymptotic bound $\OO(k^{-1})$ is therefore as small as possible.
\end{proof}

\par In the specific case of fractionally integrated noise, we may write the prediction error as:
\begin{displaymath}
\mathbb{E}\big( \big[X_{k+1}-\widetilde{X'_k}(1)\big]^2 \big) 
 = \sigma_{\e}^2+C(d) k^{-1}+\mathrm{o}\left( k^{-1}\right) 
\end{displaymath}
and we can express $C(d)$ as a function of $d$:
\begin{equation}
C(d)=\frac{\Gamma(1-2d)\Gamma(2d)}{\Gamma(-d)^2\Gamma(d)\Gamma(1+d)}\label{C}.
\end{equation}
 It is easy to prove that $C(d) \rightarrow + \infty$ as $d \rightarrow 1/2$ and we may write the following asymptotic equivalent as $d \rightarrow 1/2$:
 \begin{equation}
 C(d) \sim \frac{1}{(1-2d)\Gamma(-1/2)^2\Gamma(1/2)\Gamma(3/2)}.\label{eqcd}
 \end{equation}
As $d \rightarrow 0 $, $C(d) \rightarrow 0$ and we have the following equivalent as $d \rightarrow 0 $:
 \begin{displaymath}
 C(d) \sim d^2.
 \end{displaymath}

\begin{figure}[h]
\begin{center}
\caption{Behaviour of constant $C(d)$, $d \in [0,1/2[$, defined in \eqref{C}}
\label{figureC}
\includegraphics[width=10cm]{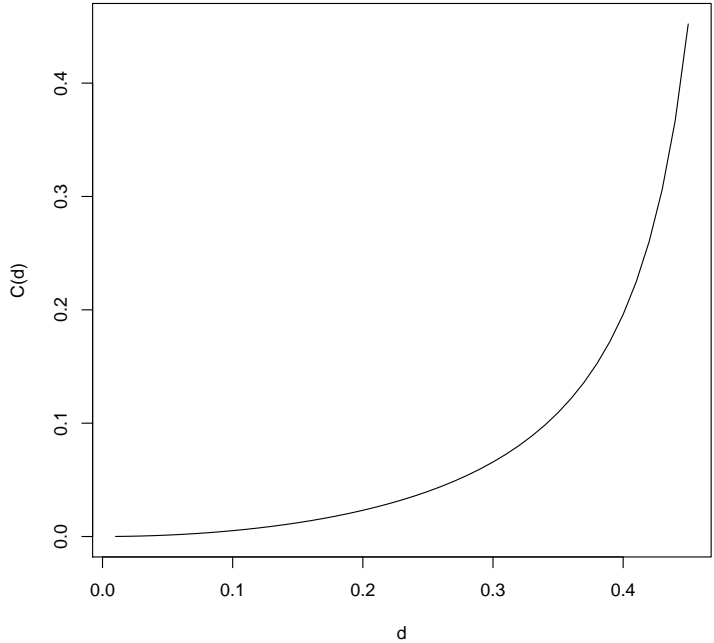}
\end{center}
\end{figure}
\par As the figure \ref{figureC} suggests and the asymptotic equivalent given in \eqref{eqcd} proves, the mean-squared error tends to $+\infty$ as $d \rightarrow 1/2$. By contrast, the constant $C(d)$ takes small values for d in a large interval of $[0,1/2[$. Although the rate of convergence has a constant order $k^{-1}$, the forecast error is bigger when $d \rightarrow 1/2$. This result is not surprising since the correlation between the random variable, which we want to predict, and the random variables, which we take equal to 0, increases when $d \rightarrow 1/2$.

\paragraph*{} Truncating to $k$ terms the series which defines the Wiener-Kolmogorov predictor amounts to using an AR($k$) model for predicting. Therefore in the following section we look for the AR($k$) which minimizes the forecast error.

\section{The Autoregressive Models Fitting Approach} \label{sectionarp}

In this section we develop a generalisation of the ``autoregressive model fitting'' approach developed by \cite{arp} in the case of fractionally integrated noise F($d$) (defined in \eqref{Fd}). We study the asymptotic properties of the forecast mean-squared error when we fit a misspecified AR($k$) model to the long-memory time series $(X_n)_{n \in \mathbb{Z}}$. 

\subsection{Rationale}

Let $\Phi$ a $k^{\textrm{th}}$ degree polynomial defined by: 
\begin{displaymath}
\Phi(z)=1-a_{1,k} z-\ldots -a_{k,k} z^k.
\end{displaymath}
We assume that $\Phi$ has no zeroes on the unit disk. We define the process $(\eta_n)_{n\in \mathbb{Z}}$ by:
\begin{displaymath}
\forall n \in \mathbb{Z}\textrm{, }\eta_n =\Phi(B) X_n
\end{displaymath} where $B$ is the backward shift operator.
Note that $(\eta_n)_{n\in \mathbb{Z}}$ is not a white noise series because $(X_n)_{n\in \mathbb{Z}}$ is a long-memory process and hence does not belong to the class of autoregressive processes. Since $\Phi$ has no root on the unit disk, $(X_n)_{n\in \mathbb{Z}}$ admits a moving-average representation as the fitted AR($k$) model in terms of $(\eta_n)_{n\in \mathbb{Z}}$:
\begin{displaymath}
X_n=\sum_{j=0}^{\infty} c(j)\eta_{n-j}.
\end{displaymath}
If $(X_n)_{n\in \mathbb{Z}}$ was an AR($k$) associated with the polynomial $\Phi$, the best next step linear predictor would be:
\begin{eqnarray*}
\widehat{X_n}(1)&=& \sum_{j=1}^{\infty} c(i) \eta_{t+1-i}\\
&=& a_{1,k}X_n +\ldots+a_{k,k} X_{n+1-k} \textrm{ si }n \geqslant k.
\end{eqnarray*}
Here $(X_n)_{n\in \mathbb{Z}}$ is a long-memory process which verifies the assumptions of Section \ref{intro}. Our goal is to derive a closed formula for the polynomial $\Phi$ which minimizes the forecast error and to estimate this error.

\subsection{Mean-Squared Error}
\label{yaj}
There exists two approaches in order to define the coefficients of the $k^{\mathrm{th}}$ degree polynomial $\Phi$: the spectral approach and the time approach. 
\par In the time approach, we choose to define the predictor as the projection mapping on to the closed span of the subset $\{X_n,\ldots,X_{n+1-k}\}$ of the Hilbert space L$^2(\Omega,\mathcal{F},\mathbb{P})$ with inner product $<X,Y>=\E(XY)$. Consequently the coefficients of $\Phi$ verify the equations, which are called the $k^{\mathrm{th}}$ order Yule-Walker equations:
\begin{equation}
\forall j \in \llbracket1,k \rrbracket, \quad \sum_{i=1}^k a_{i,k}\sigma (i-j)=\sigma(j) \label{yw}
\end{equation}

\par  The mean-squared prediction error is:
\begin{displaymath}
\mathbb{E}\big[ \big(\widehat{X_n}(1)-X_{n+1}\big)^2\big]=c(0)^2\mathbb{E}(\eta_{n+1}^2)=\mathbb{E}(\eta_{n+1}^2).
\end{displaymath}

We may write the moving average representation of $(\eta_n)_{n\in \mathbb{N}}$ in terms of $(\varepsilon_n)_{n\in \mathbb{N}}$:
\begin{eqnarray*}
\eta_n&=& \sum_{j=0}^{\infty} \sum_{k=0}^{min(j,p)}\Phi_k b(j-k) \varepsilon_{n-j}\\
&=& \sum_{j=0}^{\infty} t(j)\varepsilon_{n-j}
\end{eqnarray*} 
with
\begin{displaymath}
\forall j \in \mathbb{N}, \quad t(j)=\sum_{k=0}^{min(j,p)}\Phi_k b(j-k).
\end{displaymath}

Finally we obtain:
\begin{displaymath}
\mathbb{E}\big[ \big(\widehat{X_n}(1)-X_{n+1}\big)^2\big]=\sum_{j=0}^{\infty} t(j)^2 \sigma^2_{\varepsilon}.
\end{displaymath}
In the spectral approach, minimizing the prediction error is equivalent to minimizing a contrast between two spectral densities:
\begin{equation*}
\int_{-\pi}^{\pi}\frac{f(\lambda)}{g(\lambda,\Phi)}\dd \lambda \label{distance}
\end{equation*} 
where $f$ is the spectral density of $X_n$ and $g(.,\Phi)$ is the spectral density of the AR(p) process defined by the polynomial $\Phi$ (see for example \cite{dist}),so:
\begin{align*}
\int_{-\pi}^{\pi}\frac{f(\lambda)}{g(\lambda,\Phi)}\dd \lambda &= \int_{-\pi}^{\pi}\Big|\sum_{j=0}^{\infty}b(j)e^{-ij\lambda}\Big|^2\Big|\Phi(e^{-i\lambda })\Big|^2\dd \lambda\\
&=\int_{-\pi}^{\pi}|\sum_{j=0}^{\infty}t(j)e^{-ij\lambda}|^2\dd \lambda\\
&=2\pi\sum_{j=0}^{\infty}t(j)^2.
\end{align*}
In both approaches we need to minimize $\sum_{j=0}^{\infty}t(j)$.

\subsection{Rate of Convergence of the Error by AR($k$) Model Fitting}

In the next theorem we derive an asymptotic expression for the prediction error by fitting autoregressive models to the series:
\begin{thrm}
Assume that $(X_n)_{n \in \mathbb{Z}}$ is a long-memory process which verifies the assumptions of Section \ref{intro}. If $0<d<\frac{1}{2}$:
\begin{displaymath}
\mathbb{E}\big[ \big(\widehat{X_k}(1)-X_{k+1}\big)^2\big]-\sigma_{\e}^2=\OO(k^{-1})
\end{displaymath}
\end{thrm}
\begin{proof}
Since fitting an AR($k$) model minimizes the forecast error using $k$ observations, the error by using truncation is bigger. Since the truncation method involves an error bounded by $\OO\left(k^{-1} \right) $, we obtain:
\begin{displaymath}
\mathbb{E}\big[ \big(\widehat{X_k}(1)-X_{k+1}\big)^2\big]-\sigma_{\e}^2=\OO(k^{-1}).
\end{displaymath}
Consequently we only need to prove that this rate of convergence is attained
. This is the case for the fractionally integrated processes defined in \eqref{Fd}. We want the error made when fitting an AR($k$) model in terms of the Wiener-Kolmogorov truncation error. Note first that the variance of the white noise series is equal to:
\begin{displaymath}
\sigma_{\e}^2= \int_{-\pi}^{\pi}f(\lambda) \left| \sum_{j=0}^{+\infty}a_j e^{ij\lambda} \right|^2 \dd \lambda.
\end{displaymath}
Therefore in the case of a fractionally integrated process F($d$) we need only show that:
\begin{displaymath}
\int_{-\pi}^{\pi}f(\lambda) \left| \sum_{j=0}^{+\infty}a_j e^{ij\lambda} \right|^2 \dd \lambda-\frac{\sigma_{\e}^2}{2\pi}\int_{-\pi}^{\pi}\frac{f(\lambda)}{g(\lambda,\Phi_k)} \dd \lambda\sim C(k^{-1}).
\end{displaymath}
\begin{eqnarray*}
\int_{-\pi}^{\pi}f(\lambda) \left| \sum_{j=0}^{+\infty}a_j e^{ij\lambda} \right|^2 \dd \lambda-\frac{\sigma_{\e}^2}{2\pi}\int_{-\pi}^{\pi}\frac{f(\lambda)}{g(\lambda,\Phi_k)}\dd \lambda&=& \int_{-\pi}^{\pi}f(\lambda) \left( \left| \sum_{j=0}^{+\infty}a_j e^{ij\lambda} \right|^2- \left| \sum_{j=0}^{k}a_{j,k} e^{ij\lambda} \right|^2\right) \dd \lambda \\
&=& \sum_{j=0}^{+\infty}\sum_{l=0}^{+\infty}\left(a_ja_l-a_{j,k}a_{l,k}\right) \sigma(j-l)
\end{eqnarray*}
we set $a_{j,k}=0$ if $j>k$.
\begin{eqnarray}
&&\sum_{j=0}^{+\infty}\sum_{l=0}^{+\infty}\left(a_j a_l -a_{j,k}a_{l,k}\right) \sigma(j-l) \\ 
&=&\sum_{j=0}^{+\infty}\sum_{l=0}^{+\infty}(a_ja_l-a_{j,k}a_l)\sigma(j-l) 
+\sum_{j=0}^{+\infty}\sum_{l=0}^{+\infty}(a_{j,k}a_l-a_{j,k}a_{l,k})\sigma(j-l) \nonumber\\
&=&\sum_{j=0}^{+\infty}(a_j-a_{j,k})\sum_{l=0}^{+\infty}a_l \sigma(l-j)
+\sum_{j=0}^{k}a_{j,k}\sum_{l=0}^{+\infty}(a_l-a_{l,k})\sigma(j-l)  \quad \label{somme}
\end{eqnarray}
We first study the first term of the sum \eqref{somme}. For any  $j>0$ , we have  $\sum_{l=0}^{+\infty}a_l \sigma(l-j)=0$:
\begin{eqnarray*}
\e_n&=&\sum_{j=0}^{\infty} a_l X_{n-l}\\
X_{n-j} \e_n&=&\sum_{l=0}^{\infty} a_lX_{n-l} X_{n-j} \\
\E\left( X_{n-j} \e_n\right) &=&\sum_{l=0}^{\infty} a_l\sigma(l-j) \\
\E\left( \sum_{l=0}^{\infty} b_l \e_{n-j-l}\e_n\right)&=&\sum_{l=0}^{\infty} a_l\sigma(l-j)
\end{eqnarray*}
and we conclude that $\sum_{l=0}^{+\infty}a_l \sigma(l-j)=0$ because $(\varepsilon_n)_{n\in \mathbb{Z}}$ is an uncorrelated white noise.
We can thus rewrite the first term of \eqref{somme} like:
\begin{eqnarray}
\sum_{j=0}^{+\infty}(a_j-a_{j,k})\sum_{l=0}^{+\infty}a_l \sigma(l-j)&=&(a_0-a_{0,k})\sum_{l=0}^{+\infty}a_l \sigma(l) \nonumber \\
&=&0 \nonumber 
\end{eqnarray}
since $a_0=a_{0,k}=1$ according to definition. 
Next we study the second term of the sum \eqref{somme}:
\begin{eqnarray}
\sum_{j=0}^{k}a_{j,k}\sum_{l=0}^{+\infty}(a_l-a_{l,k})\sigma(j-l) \nonumber .
\end{eqnarray}
And we obtain that:
\begin{eqnarray}
\sum_{j=0}^{k}a_{j,k}\sum_{l=0}^{+\infty}(a_l-a_{l,k})\sigma(j-l)&=& \sum_{j=1}^{k}(a_{j,k}-a_j)\sum_{l=1}^{k}(a_l-a_{l,k})\sigma(j-l) \nonumber \\
&&+\sum_{j=1}^{k}(a_{j,k}-a_j)\sum_{l=k+1}^{+\infty}a_l\sigma(j-l) \label{milieu1}\\
&&+\sum_{j=0}^{k}a_j\sum_{l=1}^{k}(a_l-a_{l,k})\sigma(j-l) \label{milieu2}\\
&&+\sum_{j=0}^{k}a_j\sum_{l=k+1}^{+\infty}a_l\sigma(j-l) \nonumber
\end{eqnarray}
Similarly we rewrite the term \eqref{milieu1} using the Yule-Walker equations:
\begin{displaymath}
\sum_{j=1}^{k}(a_{j,k}-a_j)\sum_{l=k+1}^{+\infty}a_l\sigma(j-l)=-\sum_{j=1}^{k}(a_{j,k}-a_j)\sum_{l=0}^{k}a_l\sigma(j-l)
\end{displaymath} 
We then remark that this is equal to \eqref{milieu2}.
Hence it follows that:
\begin{eqnarray}
\sum_{j=0}^{k}a_{j,k}\sum_{l=0}^{+\infty}(a_l-a_{l,k})\sigma(j-l)&=&\sum_{j=1}^{k}(a_{j,k}-a_j)\sum_{l=1}^{k}(a_l-a_{l,k})\sigma(j-l) \nonumber \\
&&+2\sum_{j=1}^{k}(a_{j,k}-a_j) \sum_{l=k+1}^{+\infty}a_l\sigma(j-l) \nonumber \\
&&+\sum_{j=0}^{k}a_j\sum_{l=k+1}^{+\infty}a_l\sigma(j-l) \label{commeavant}
\end{eqnarray}
On a similar way we can rewrite the third term of the sum \eqref{commeavant} using Fubini Theorem:
\begin{displaymath}
\sum_{j=0}^{k}a_j\sum_{l=k+1}^{+\infty}a_l\sigma(j-l)=-\sum_{j=k+1}^{+\infty}\sum_{l=k+1}^{+\infty}a_ja_l \sigma(j-l).
\end{displaymath}
This third term is therefore equal to the forecast error in the method of prediction by truncation.

  In order to compare the prediction error by truncating the Wiener-Kolmogorov predictor and by fitting an autoregressive model to a fractionally integrated process F($d$), we need the sign of all the components of the sum \eqref{commeavant}. For a fractionally integrated noise, we know the explicit formula for $a_j$ and $\sigma(j)$:
\begin{displaymath}
\forall j>0,\quad a_j=\frac{\Gamma(j-d)}{\Gamma(j+1)\Gamma(-d)}<0 \textrm{ and } \forall j\geq 0,\quad \sigma(j)=\frac{(-1)^j\Gamma(1-2d)}{\Gamma(j-d+1)\Gamma(1-j-d)} \sigma_\e^2>0.
\end{displaymath}
In order to get the sign of $a_{j,k}-a_j$ we use the explicit formula given in \cite{coefftronccourtemem} and we easily obtain that $a_{j,k}-a_j$ is negative for all $j \in \llbracket 1, k \rrbracket$.
\begin{eqnarray*}
a_j-a_{j,k}&=&\frac{\Gamma(j-d)}{\Gamma(j+1)\Gamma(-d)}-\frac{\Gamma(k+1)\Gamma(j-d) \Gamma(k-d-j+1)}{\Gamma(k-j+1)\Gamma(j+1)\Gamma(-d)\Gamma(k-d+1)}\\
&=&-a_j\left(-1+\frac{\Gamma(k+1)\Gamma(k-d-j+1)}{\Gamma(k-j+1)\Gamma(k-d+1)} \right) \\
&=&-a_j\left(\frac{k...(k-j+1)}{(k-d)...(k-d-j+1)}-1\right) \\
&>&0
\end{eqnarray*}
since $\forall j \in \mathbb{N}^*$ $a_j<0$.
 To give an asymptotic equivalent for the prediction error, we use the sum given in \eqref{commeavant}. We have the sign of the three terms: the first is negative, the second is positive and the last is negative. Moreover the third is equal to the forecast error by truncation and we have proved that this asymptotic equivalent has order $\OO(k^{-1})$. The  prediction error by fitting an autoregressive model converges faster to 0 than the error by truncation only if the second term is equivalent to $Ck^{-1}$, with $C$ constant. Consequently, we search for a bound for $a_j-a_{j,k}$ given the explicit formula for these coefficients (see for example \cite{coefftronccourtemem}):
\begin{eqnarray*}
a_j-a_{j,k}&=&\frac{\Gamma(j-d)}{\Gamma(j+1)\Gamma(-d)}-\frac{\Gamma(k+1)\Gamma(j-d) \Gamma(k-d-j+1)}{\Gamma(k-j+1)\Gamma(j+1)\Gamma(-d)\Gamma(k-d+1)}\\
&=&-a_j\left(-1+\frac{\Gamma(k+1)\Gamma(k-d-j+1)}{\Gamma(k-j+1)\Gamma(k-d+1)} \right) \\
&=&-a_j\left(\frac{k...(k-j+1)}{(k-d)...(k-d-j+1)}-1\right) \\
&=&-a_j\left(\prod_{m=0}^{j-1} \left( \frac{1-\frac{l}{k}}{1-\frac{l+d}{k}}\right) -1\right) \\
&=&-a_j\left(\prod_{m=0}^{j-1}\left( 1+\frac{\frac{d}{k}}{1-\frac{d+l}{k}}\right) -1\right) .
\end{eqnarray*}
Then we use the following inequality:
\begin{displaymath}
\forall x \in \mathbb{R}, \quad 1+x\leq \exp(x)
\end{displaymath}
which gives us:
\begin{eqnarray*}
a_j-a_{j,k}&\leq&-a_j\left(\exp\left(\sum_{m=0}^{j-1} \frac{\frac{d}{k}}{1-\frac{d+l}{k}}\right) -1\right) \\
&\leq&-a_j\left(\exp \left(d\sum_{m=0}^{j-1}\frac{1}{k-d-l} \right) -1\right)\\
&\leq&-a_j\exp \left(d\sum_{m=0}^{j-1}\frac{1}{k-d-l} \right) 
\end{eqnarray*}
According to the previous inequality, we have:
\begin{eqnarray*}
\sum_{j=1}^{k}(a_j-a_{j,k})\sum_{l=k+1}^{+\infty}-a_l \sigma(j-l)&=&\sum_{j=1}^{k-1}(a_j-a_{j,k})\sum_{l=k+1}^{+\infty}-a_l \sigma(j-l) \\
&&+(a_k-a_{k,k})\sum_{l=k+1}^{+\infty}-a_l \sigma(k-l) \\
&\leq & \sum_{j=1}^{k-1}-a_j\exp \left(d\sum_{m=0}^{j-1}\frac{1}{k-d-m}\right)  \sum_{l=k+1}^{+\infty}-a_l \sigma(j-l) \\
&&+(-a_k)\exp \left(d\sum_{m=0}^{k-1}\frac{1}{k-d-m} \right)\sum_{l=k+1}^{+\infty}-a_l \sigma(k-l) \\
&\leq & \sum_{j=1}^{k-1}-a_j\exp \left(d\int_0^{j}\frac{1}{k-d-m} \dd m\right)\sum_{l=k+1}^{+\infty}-a_l \sigma(j-l) \\
&&+(-a_k) k^{\frac{3}{2}d}\sum_{l=k+1}^{+\infty}-a_l \sigma(k-l) 
\end{eqnarray*}
As the function $x \mapsto \frac{1}{k-d-x}$ is increasing, we use the Integral Test Theorem. The inequality on the second term follows from:
\begin{eqnarray*}
\sum_{m=0}^{k-1}\frac{1}{k-d-m} &\sim& \ln(k)\\
&\leq& \frac{3}{2}\ln(k) 
\end{eqnarray*}for $k$ large enough.
Therefore there exists $K$ such that for all $k \geq K$:
\begin{eqnarray*}
\sum_{j=1}^{k}(a_j-a_{j,k})\sum_{l=k+1}^{+\infty}-a_l \sigma(j-l)&\leq & \sum_{j=1}^{k-1}-a_j\exp \left(d \ln\left( \frac{k-d}{k-d-j}\right) \right) \sum_{l=k+1}^{+\infty}-a_l \sigma(j-l)\\
&&+(-a_k) k^{\frac{3}{2}d}\sum_{l=k+1}^{+\infty}-a_l \sigma(0)\\
&\leq &C (k-d)^d\sum_{j=1}^{k-1}j^{-d-1}(k-d-j)^{-d}\sum_{l=k+1}^{+\infty}l^{-d-1}(l-j)^{2d-1}\\
&&+Ck^{-d-1}k^{\frac{3}{2}d}k^{-d} \\
&\leq &\frac{C}{(k-d)^{2}}\int_{1/(k-d)}^1 j^{-d-1}(1-j)^{-d}\int_{1}^{+\infty}l^{-d-1}(l-1)^{2d-1} \dd l \dd j\\
&&+Ck^{-\frac{1}{2}d-1}\\
&\leq &C'(k-d)^{-2+d}+Ck^{-\frac{1}{2}d-1}
\end{eqnarray*}
and so the positive term has a smaller asymptotic order than the forecast error made by truncating. Therefore we have proved that in the particular case of F($d$) processes, the two prediction errors are equivalent to $Ck^{-1}$ with $C$ constant. 
\end{proof}

\begin{figure}[h]
\begin{center}
\caption{Ratio $r(k)$, $d \in ]0,1/2[$ defined in \eqref{differreur}}
\label{graphefigerreur}
\includegraphics[width=10cm]{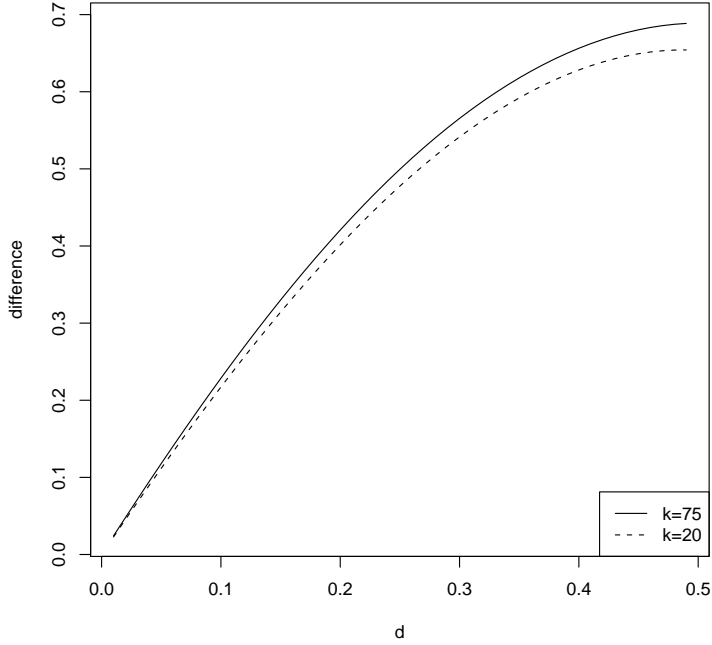}
\end{center}
\end{figure}

The two approaches to next-step prediction, by truncation to $k$ terms or by fitting an autoregressive model AR($k$) have consequently a prediction error with the same rate of convergence $k^{-1}$. So it is interesting to study how the second approach improves the prediction. The following quotient:
\begin{equation}
r(k):=\frac{\sum_{j=1}^{k}(a_{j,k}-a_j)\sum_{l=1}^{k}(a_l-a_{l,k})\sigma(j-l) 
+2\sum_{j=1}^{k}(a_{j,k}-a_j) \sum_{l=k+1}^{+\infty}a_l\sigma(j-l)}{\sum_{j=0}^{k}a_j\sum_{l=k+1}^{+\infty}a_l\sigma(j-l)} \label{differreur}
\end{equation}
is the ratio of the difference between the two prediction errors and the prediction error by truncating in the particular case of a fractionally integrated noise F($d$). The figure \ref{graphefigerreur} shows that the prediction by truncation incurs a larger performance loss when $d \rightarrow 1/2$. The improvement reaches 50 per cent when $d>0.3$ and $k>20$.

\par After obtaining asymptotic equivalent for next step predictor, we will generalize the two methods of $h$-step prediction and aim to obtain their asymptotic behaviour as $k \rightarrow + \infty$ but also as $h\rightarrow + \infty$.

\section{The h-Step Predictors}
\label{prevh}
Since we assume that the process $(X_n)_{n \in \mathbb{Z}}$ has an autoregressive representation \eqref{defaj} and moving average representation \eqref{defbj}, the linear least-squares predictor, $\widetilde{X}_{k+h}$, of $X_{k+h}$ based on the infinite past $(X_j,j\leq k)$ is given by:
\begin{displaymath}
\widetilde{X_k}(h)= - \sum_{j=1}^{+\infty} a_j \widetilde{X_k}(h-j)= \sum_{j=h}^{+\infty} b_j \e_{k+h-j}
\end{displaymath}(see for example Theorem 5.5.1 of \cite{bd}).
The corresponding mean squared error of prediction is:
\begin{displaymath}
\E\left[ \left(\widetilde{X_k}(h) -X_{k+h}\right)^2 \right] = \sigma_\e^2 \sum_{j=0}^{h-1}b_j^2.
\end{displaymath}
As the prediction step $h$ tends to infinity, The mean-squared prediction error converges to $\sigma_\e^2 \sum_{j=0}^{+\infty}$, which is the the variance of the process $(X_n)_{n \in \mathbb{Z}}$. But if the mean-squared prediction error is equal to $\sigma(0)$, we have no more interest in the prediction method since its error is equal to the error of predicting the future by 0.
Remark that the mean-squared error increases more slowly to $\sigma(0)$ in the long-memory case than in the short-memory case since the sequence $b_j$ decays more slowly to 0. More precisely in the case of a long-memory process, if we assume that:
\begin{displaymath}
b_j \underset{+\infty}\sim j^{d-1} L(j)
 \end{displaymath}
 where $L$ is a slowly varying function, we can express the asymptotic behaviour of the prediction error. As $j \mapsto L^2(j)$ is also a slowly varying function according to the definition of \cite{zygmund}, $b_j^2=j^{2d-2} L^2(j)$ is ultimately decreasing. The rest of the series and the integral are then equivalent and we may write:
\begin{eqnarray}
\sigma(0)-\E\left[ \left( \widetilde{X_k}(h)-X_{k+h}\right)^2 \right] &=& \sigma_\e^2 \sum_{j=h}^{+\infty} b_j^2 \nonumber \\
& \sim &\sum_{j=h}^{+\infty}j^{2d-2} L^2(j)\nonumber\\
& \sim & \int_{h}^{+\infty} j^{2d-2} L^2(j) \dd j \nonumber
\end{eqnarray}
According to Proposition 1.5.10 of \cite{bingham}:
\begin{eqnarray}
\sigma(0)-\E\left[ \left( \widetilde{X_k}(h)-X_{k+h}\right)^2 \right]& \sim & \int_{h}^{+\infty} j^{2d-2} L^2(j) \dd j\nonumber\\
& \underset{h \rightarrow +\infty}\sim & \frac{1}{1-2d}h^{2d-1}L^2(h)\label{inf}
\end{eqnarray}
In the case of a long-memory process with parameter $d$ which verifies $b_j \sim j^{d-1}L(j)$, the convergence of the mean-squared error to $\sigma(0)$ is slow as $h$ tends to infinity.
On the contrary, for a moving average process of order $q$, the sequence $\sigma(0)-\E\left[ \left( \widetilde{X_k}(h)-X_{k+h}\right)^2 \right] $ is constant and equal to 0 as soon as $h>q$. More generally, we can study the case of an ARMA process, which canonical representation is given by:
\begin{displaymath}
\Phi(X_t)=\Theta(\e_t)
\end{displaymath}
where $\Phi$ and $\Theta$ are two coprime polynomials with coefficients of degree 0 are equal to 1 and $\e_t$ is a white noise. $\Phi$ has no root in the unit disk $|z|\leq1$ and $\Theta$ has no root in the open disk $|z|<1$.
$b_j$ is bounded by:
\begin{displaymath}
|b_j |\leq  C j^{m-1} \rho^{-j} 
\end{displaymath}
where $\rho$ is the smallest absolute value of the roots of $\Phi$ and $m$ the multiplicity of the corresponding root (see for example \cite{bd} p92). Thus the mean-squared prediction error is bounded by:
\begin{eqnarray*}
\sigma(0)-\E\left[ \left( \widetilde{X_k}(h)-X_{k+h}\right)^2 \right] &=& \sigma_\e^2 \sum_{j=h}^{+\infty} b_j^2 \\
&\leq & \sigma_\e^2 C^2 \sum_{j=h}^{+\infty}j^{2m-2}\rho^{-2j}\\
&\leq & \sigma_\e^2 C^2\sum_{j=h}^{+\infty}j^{2m-2} \exp\left( -2j\log({\rho})\right)  \\
&\leq &\sigma_\e^2 C^2\int_{h}^{+\infty}j^{2m-2} \exp\left( -2j\log({\rho})\right) \dd j 
\end{eqnarray*}
By using the substitution $t=2\log({\rho})j$ ,
\begin{eqnarray*}
 \sigma(0)-\E\left[ \left(\widetilde{X_k}(h)-X_{k+h}\right)^2 \right] &\leq &\sigma_\e^2 C^2\left( 2\log({\rho})\right)^{1-2m}\int_{2\log({\rho})h}^{+\infty}t^{2m-2} \exp\left( t\right) \dd t \\
&\leq &\sigma_\e^2 C^2\left( 2\log({\rho})\right)^{1-2m} \Gamma(2m-1,2\log({\rho})h)
\end{eqnarray*}
where $\Gamma(.,.)$ is the incomplete Gamma function defined in equation 6.5.3 of \cite{abramo}. We know an equivalent of this function:
\begin{displaymath}
  \Gamma(2m-1,2\log({\rho})h) \underset{h \rightarrow +\infty}\sim\left(2\log({\rho})h \right) ^{2m-2} \exp\left(2\log({\rho})h \right)
\end{displaymath}
We conclude that the rate of convergence is exponential. 
The mean-squared prediction error goes faster to $\sigma(0)$ when the predicting process is ARMA than when the process is a long-memory process.\\
The $h$-step prediction is then more interesting for the long-memory process than for the short-memory process, having observed the infinite past. We consider the truncating effect next.

\subsection{Truncated Wiener-Kolmogorov predictor}
In practice, we only observe a finite number of samples. We assume now that we only know $k$ observations $(X_1,\ldots,X_k)$.
We then define the $h$-step truncated Wiener-Kolmogorov of order $k$ as:
\begin{equation}
\widetilde{X'_k}(h)= -\sum_{j=1}^{h-1}a_j \widetilde{X'_k}(h-j) -\sum_{l=1}^k a_{h-1+j} X_{k+1-j}\label{wkth}
\end{equation}

We now describe the asymptotic behaviour of the mean-squared error of the predictor \eqref{wkth}. First we write the difference between the predicting random variable and its predictor:
\begin{eqnarray*}
\widetilde{X'_k}(h)-X_{k+h}&=&-\sum_{j=1}^{h-1}a_j \widetilde{X'_k}(h-j) -\sum_{l=1}^k a_{h-1+j} X_{k+1-j}-\e_{k+h} +\sum_{j=1}^{+\infty}a_jX_{k+h-j}
\\&=&-\e_{k+h}+\sum_{j=1}^{h-1}a_j \left( X_{k+h-j}-\widetilde{X'_k}(h-j)\right)+ \sum_{j=1}^k a_{h-1+j}\left(X_{k+1-j}- X_{k+1-j}\right)\\
&& +\sum_{j=k+1}^{+\infty} a_{h-1+j}X_{k+1-j}\\
&=&-\e_{k+h}+\sum_{j=1}^{h-1}a_j \left( X_{k+h-j}-\widetilde{X'_k}(h-j)\right)+\sum_{j=k+1}^{+\infty} a_{h-1+j}X_{k+1-j}
\end{eqnarray*}
We will use the process of induction on $h$ to show that 
\begin{eqnarray}
\widetilde{X'_k}(h)-X_{k+h}&=&-\sum_{l=0}^{h-1} \left( \sum_{j_1+j_2+\ldots+j_h=l} (-1)^{card(\{j,j\neq0\})} a_{j_1} a_{j_2} \ldots a_{j_h} \right) \e_{k+h-l} \nonumber\\
&&+\sum_{j=k+1}^{+\infty} \left(\sum_{i_1+i_2+\ldots+i_{h}=h-1}(-1)^{card(\{i_l,i_l\neq0,l>1\})} a_{j+i_1} a_{i_2} \ldots a_{i_h} \right)X_{k+1-j}. \nonumber
\end{eqnarray} 
For $h=2$, we have for example
\begin{displaymath}
\widetilde{X'_k}(2)-X_{k+2}=-(a_0\e_{k+2}-a_1\e_{k+1})+\sum_{j=k+1}^{+\infty}(-a_1a_j+a_{j+1})X_{k+1-j}.
\end{displaymath}
Let $A(z)$ and $B(z)$ denote $A(z)=1+ \sum_{j=1}^{+\infty}a_j z^j$ and $B(z)=1+ \sum_{j=1}^{+\infty}b_j z^j$. Since we have $A(z)=B(z)^{-1}$, we obtain the following conditions on the coefficients:
\begin{eqnarray*}
b_1&=&-a_1\\
b_2&=& -a_2+a_1^2\\
b_3&=& -a_3+2a_1a_2-a_1^3\\
\ldots
\end{eqnarray*}
So we obtain:
\begin{equation}
\widetilde{X'_k}(h)-X_{k+h}=-\sum_{l=0}^{h-1} b_l\e_{k+h-l} +\sum_{j=k+1}^{+\infty} \sum_{m=0}^{h-1} a_{j+m}b_{h-1-m}X_{k+1-j} \label{suminnovtronc}.
\end{equation}
Since the process $(\e_n)_{n \in \mathbb{Z}}$ is uncorrelated and then the two terms of the sum \eqref{suminnovtronc} are orthogonal, we can rewrite the mean-squared error:
\begin{eqnarray}
\E\left[ \widetilde{X'_k}(h)-X_{k+h}\right]^2
&=& \sum_{l=0}^{h-1} b_l^2 \sigma_\e^2 \label{methode}\\&&+\E\left[\sum_{j=k+1}^{+\infty} \left( \sum_{m=0}^{h-1} a_{j+h-1-m}b_{m} \right)X_{k+1-j}\right]^2.\label{tronc}
\end{eqnarray}
The first part of the error \eqref{methode} is due to the prediction method and the second \eqref{tronc} due to the truncating of the predictor. 
We now approximate the error term \eqref{tronc} by using \eqref{aj} and \eqref{bj}. We obtain the following upper bound:
\begin{eqnarray}
\forall \delta >0,\:\left\vert \sum_{m=0}^{h-1} a_{j+h-1-m}b_{m}\right\vert&\leq& \sum_{m=1}^{h-1}\vert a_{j+h-1-m}b_{m}\vert +\vert b_0 a_{j+h-1} \vert\nonumber\\
& \leq& C_1C_2\int_{0}^{h}(j+h-1-l)^{-d-1+\delta}l^{d-1+\delta} \dd l +C_1 (j+h)^{-d-1}\nonumber\\
& \leq&C_1C_2 h^{-1+2\delta}\int_{0}^{1} \left(\frac{j}{h}+1-l\right)^{-d-1+\delta}l^{d-1+\delta}\dd l+ C_1 (j+h)^{-d-1}\nonumber\\
& \leq&C_1C_2 h^{-1+2\delta} j^{-d-1+\delta}\int_{0}^{1}\left(\frac{1}{h}+\frac{1-l}{j}\right)^{-d-1+\delta}l^{d-1+\delta} \dd l+C_1 (j+h)^{-d-1}\nonumber\\
& \leq&C_1C_2h^{d+2\delta}j^{-d-1+\delta}\int_{0}^{1}l^{d-1+\delta} \dd l+C_1 (j+h)^{-d-1}\nonumber\\
\left\vert \sum_{m=0}^h a_{j+h-m}b_{m}\right\vert& \leq&C_1C_2\frac{h^{d+2\delta}}{d}j^{-d-1+\delta} \label{cjh}
\end{eqnarray}
This bound is in fact an asymptotic equivalent for the fractionally integrated process F($d$) because, in that case, the sequences $a_j$ and $b_j$ have a constant signs. Using Proposition \ref{erreurtronque} for the one-step prediction and we have:
\begin{prpstn}Let $(X_n)_{n \in \mathbb{Z}}$ be a linear stationary process defined by \eqref{defbj}, \eqref{defaj} and possessing the features \eqref{aj} and \eqref{bj}. We can approximate the mean-squared prediction error of $\widetilde{X'_k}(1)$ by:
\begin{equation}
\forall \delta>0, \:\E\left[ \widetilde{X'_k}(h)-X_{k+h}\right]^2=\sum_{l=0}^{h-1} b_l^2 \sigma_\e^2+\OO\left( h^{2d+\delta} k^{-1+\delta}\right).\label{tronch}
\end{equation}
\end{prpstn}
 Having $k$ observations, we search for the step $h$ for which the variance of the predictor has for upper bound $\sigma(0)$. Then the prediction error have for asymptotic bound  $\OO\left(h^{2d}k^{-1} \right)$. We want to choose $h$ to have the prediction error negligible with respect to the information given by the linear least-squares predictor given the infinite past (see \eqref{inf}) and we obtain:
\begin{displaymath}
h^{2d}k^{-1}=\mathrm{o}(h^{2d-1})
\end{displaymath}
and then $h=\mathrm{o}(k)$. With the truncated Wiener-Kolmogorov predictor, it is interesting to compute the $h$-step predictor if we have $k$ observations $h=\mathrm{o}(k)$.

\subsection{The k-th Order Linear Least-Squares Predictor}
 For next step predictor, when we fitted an autoregressive process, we search the linear least-squares predictor knowing the finite past $(X_1,\ldots,X_k)$ and the predictor is then the projection of the random variable onto the past. Let $\widehat{X_k}(h)$ denote the projection of $X_{k+h}$ onto the span of $(X_1,\ldots,X_k)$. $\widehat{X_k}(h)$ verifies the recurrence relationship $$\widehat{X_k}(h) =-\sum_{j=1}^k a_{j,k} \widehat{X_k}(h-j)$$ where $\widehat{X_k}(h-j)$ is the direct linear least-squares predictor of $X_{k+h-j}$ based on the finite past $(X_1, \ldots,X_k)$. By induction, we obtain the predictor as a function of $(X_1, \ldots,X_k)$:
For next step prediction by fitting an autoregressive process, the best linear least-squares predictor knowing the finite past is a projection of the random variable $X_{k+1}$ onto the past. 
\begin{displaymath}
\widehat{X_k}(h) = -\sum_{j=1}^k c_{j,k} X_{k+1-j}.
\end{displaymath}
Since $\widehat{X_k}(h) $ is the projection of $X_{k+h}$ onto $(X_1, \ldots,X_k)$ in $\mathrm{L}^2$, the vector $(c_{j,k})_{1\leq j\leq k}$ minimizes the mean-squared error:
\begin{displaymath}
\E\left[ \widehat{X_k}(h)-X_{k+h}\right]^2=\int_{-\pi}^{\pi}f(\lambda) \left \vert \exp(\ii \lambda(h-1)) +\sum_{j=1}^k c_{j,k} \exp(-\ii \lambda j) \right\vert ^2\dd \lambda 
\end{displaymath}
The vector $(c_{j,k})_{1\leq j\leq k}$ is a solution of the equation:
\begin{displaymath}
\nabla_c \,\E\left[ \widehat{X_k}(h)-X_{k+h}\right]^2=0
\end{displaymath}
where $\nabla_c$ is the gradient. The vector $(c_{j,k})_{1\leq j\leq k}$ is then equal to:
\begin{equation}
(c_{j,k})_{1\leq j\leq k} = -\Sigma_k^{-1} (\sigma_{h-1+j})_{1\leq j\leq k}.\label{ywh}
\end{equation}

The corresponding mean squared error of prediction is given by:
\begin{eqnarray*}
\E\left[ \widehat{X_k}(h)-X_{k+h}\right]^2&=&\int_{-\pi}^{\pi}f(\lambda) \left \vert \exp(\ii \lambda(h-1)) +\sum_{j=1}^k c_{j,k} \exp(-\ii \lambda j) \right\vert ^2\dd \lambda \\
&=& \sigma(0)+2\sum_{j=1}^k c_{j,k} \sigma(h-1+j) +\sum_{j,l=1}^k c_{j,k} c_{l,k} \sigma(j-l) \\
&=& \sigma(0)+2\;{}^t(c_{j,k})_{1\leq j\leq k}(\sigma_{h-1+j})_{1\leq j\leq k}+{}^t(c_{j,k})_{1\leq j\leq k}\Sigma_k(c_{j,k})_{1\leq j\leq k}\\
&=& \sigma(0)- {}^t(\sigma_{h-1+j})_{1\leq j\leq k}\Sigma_k^{-1} (\sigma_{h-1+j})_{1\leq j\leq k}
\end{eqnarray*}
The matrix $\Sigma_k^{-1}$ is symmetric positive definite and the prediction error of this method is always lower than $\sigma(0)$.

\begin{figure}[h]
\begin{center}
\caption{Mean-squared error of $\widetilde{X_k}(h)$ (MMSE), $\widetilde{X'_k}(h)$ (TPMSE) and $\widehat{X_k}(h)$ (LLSPE) for $d=0.4$ and $k=80$}
\label{graphprevh}
\includegraphics[width=10cm,angle=-90]{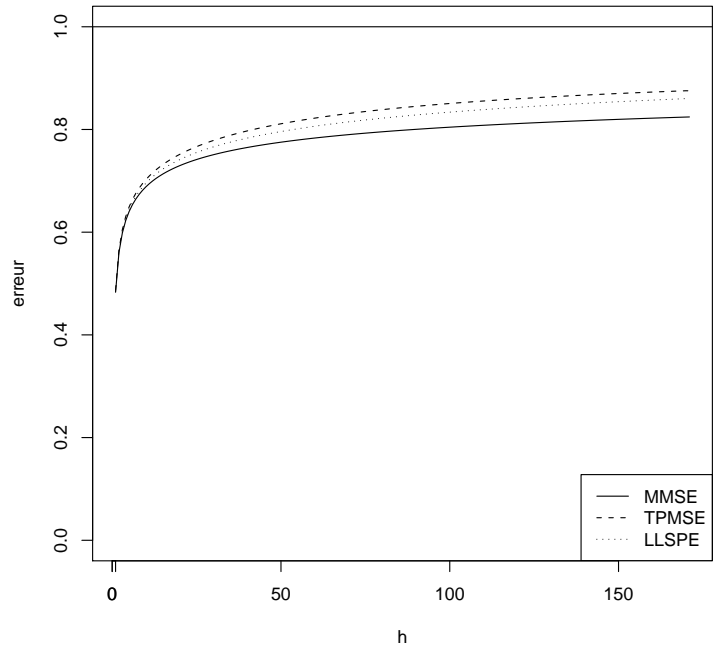}
\end{center}
\end{figure}

As $\widehat{X_k}(h)$ is the projection of $X_{k+h}$ onto $(X_1,\ldots,X_k)$,  the mean-squared prediction error is also lower than the prediction error of the truncated Wiener-Kolmogorov predictor (see figure \ref{graphprevh}). The mean-squared error of prediction due to the projection onto the span of $(X_1,\ldots, X_k)$ tends at least as fast to zero as the mean-squared due to truncation of the least-squares predictor. For one-step predictor,we have shown that the two methods can have the same rate of convergence.
\clearpage

\bibliographystyle{apalike} 
\bibliography{biblioprev}

\end{document}